# The oldest example of compound interest in Sumer:
# Seventh power of four-thirds

Kazuo MUROI

## §1. The origin of compound interest

In many modern countries – except Islamic – compound interest is used to compute interest on a loan, and plays an important role in banking and finance. The compounding of interest, which is often described as "interest on interest", is occasionally considered to be one of the greatest inventions by human beings, although a high rate of compound interest can cause serious social problems.

It is generally agreed that the origin of compound interest can be traced back to the Old Babylonian period (ca. 2000–1600 BCE), because we know that the Babylonians called compound interest *ṣibāt ṣibtim* "interest on interest" in Akkadian, and even solved mathematical problems on it. See §3 below. There is, however, a possibility that its origin goes back even further – to the Pre-Sargonic period (ca. 2600–2350 BCE), for we also know the Sumerian words maš or máš "interest" and $ur_5$-(ra) "interest-bearing loan" were used in those days.

In the present paper I shall offer a strong piece of evidence that compound interest *began at Sumer*, clarifying two passages of the Enmetena Foundation Cone[1] for which no reasonable interpretation has yet been given.

## §2. Compound interest on barley

Around 2400 BCE, the ruler of the Sumerian city of Lagash, Enmetena, had the details of the battle against the neighboring city, Umma, recorded on the foundation clay cone of his new temple. According to the following passage inscribed on the Cone,[2] Lagash had lent a large amount of barley at interest to Umma, but Umma did not repay

the loan:

| | |
|---|---|
| 21. 1 gur$_7$-am$_6$ | 1 gur$_7$ (of barley) |
| 22. lú Umma$^{ki}$-ke$_4$ | the man of Umma got an |
| 23. ur$_5$-šè ì-kú | interest-bearing loan (from Lagash). |
| 24. kud-rá ba-ús | The decision followed. |
| 25. 4 šar'u-gal gur$_7$ | (The principal and interest) became |
| 26. ba-ku$_4$ | 40,0,0,0 (sìla, whose order is) gur$_7$. |

In order to facilitate the understanding of the numerical relation between lines 21 and 25, several explanations for the capacity units, numerals and the interest rate are given below.

(1) The capacity unit gur$_7$

It is the largest capacity unit, and is defined by the basic capacity unit sìla (≈ 1 litre) as:

1 gur$_7$ = 5,20,0,0 sìla (= 1152000 sìla).

This relation is clearly confirmed by the mathematical exercises of those days, that is, by TSŠ 50 and 671.[3]

(2) The numeral šar'u-gal

Although we do not know the precise pronunciation of Sumerian numerals, we definitely know that the following numerals for large numbers were often used in administrative texts.

šár = 1,0,0 (= 3600)

šar'u (ŠÁR × U) = 10,0,0 (= 36000)

šár-gal = 1,0,0,0 (= 216000), where the adjective gal means "big".[4]

We can, therefore, interpret 4 šar'u-gal to be;

40,0,0,0 (sìla) or 8640000 (sìla).

The sign gur₇ at the end of line 25 doesn't mean 40,0,0,0 gur₇, which is too large, but that 40,0,0,0 (sìla) amounts to some of the barley in gur₇-unit, that is, 40,0,0,0 (sìla =7.5) gur₇. In mathematical texts of the Old Babylonian period we find the same custom of writing units of length in this way, which may seem odd to modern readers:

30 gi "0;30 (nindan, that is, 1) gi"

where the basic unit of length is nindan (≈ 6m), which is also omitted in the text, and

1 gi = 0;30 nindan.[5]

Since neither the Babylonians nor the Sumerians invented the number zero nor "the sexagesimal point", they occasionally resorted to this clumsy writing method in order to clarify the absolute value of a number or to read a number correctly.

(3) The rate of interest

In Mesopotamia, the rate of interest on barley was usually 33$^{1}/_{3}$% as the codes of the Old Babylonian period prescribe:

[šu]m-ma da[m-gàr] kù$^{sic}$-babbar$^{sic}$-am / a-na ur₅-ra [id-d]i-in / a-na 1 gur-e [1 (PI) 4 (bán)] še máš / i-le-[e]q-qé

"If a merchant gives silver$^{sic}$ (mistake for barley) on interest-bearing loan, he will take [1 (PI) 4 (bán) (= 60 + 40 sìla)] of barley per 1 gur (= 300 sìla) as interest."[6]

The interest rate is therefore $\dfrac{(60+40)}{300}$ = 33$^{1}/_{3}$%.

Now, we are ready to analyse the data given in lines 21 and 25. Since the 8640000 sìla in line 25 is seven and a half times the 1152000 sìla in line 21, we can assume that the former was the result of compound interest applied to the latter. So, how many years did the total take to become 7.5 times the principal? The following equation may hold in our case and it gives us a definite answer to the question above.

$$\left(1+\frac{1}{3}\right)^x = 7.5$$

$$x = \frac{1+\log 3 - 2\log 2}{2\log 2 - \log 3} \approx \frac{0.875}{0.125} = 7 \text{ years},$$

where the approximations $\log 2 \approx 0.301$ and $\log 3 \approx 0.477$ are used.

We have calculated that the Sumerian scribe wished to calculate the total due for repayment on the principal of 5,20,0,0 sìla at 33$^1$/$_3$% over seven years, i.e.

$$\left(1+\frac{1}{3}\right)^7 \cdot 5{,}20{,}0{,}0$$

It is probable that he used the approximation:

$(4/3)^7 = 7;29,29,32,50,22,13,20 \approx 7;30,$

and multiplied the principal by 7;30, that is,

7;30 · 5,20,0,0 = 40,0,0,0.

In addition, I would like to point out an important aspect of the number 7 in Sumerian superstitions. As I have already clarified elsewhere, the Sumerians considered the number 7, among the natural numbers, to be special because they believed that particular events would happen at the seventh attempt.[7] So Lagash must have waited for seven years for the repayment of barley by Umma and then have claimed the barley and its interest back.

However, the event which actually happened at the seventh year was not the repayment of barley, but an attack against Lagash by Ur-lumma, a ruler of Umma. Having repulsed the army of Umma, Enmetena finally put Ur-lumma to death in his own city of Umma. It seems that the next ruler of Umma, Il, obeyed Enmetena, for he cleared all Umma's debt to Lagash after a year, as another passage of the Cone clearly shows:[8]

| | |
|---|---|
| 11. še Lagash$^{ki}$ 10 gur$_7$-am$_6$ | 10 gur$_7$ of Lagash's barley |
| 12. ì-su | $(7.5 \cdot (4/3) = 10)$ he repaid. |

Thus the Sumerian scribe correctly calculated the compound interest on the loan of barley, whereas modern scholars, strangely enough, have not been able to understand his calculation since the end of the 19th century.

### §3. Compound interest in Babylonian mathematics

Having taken the custom of compound interest invented by the Sumerians, the Babylonians made use of it even in mathematics. Let us quote three characteristic equations of compound interest presented in Old Babylonian mathematical problems.[9] In these three problems the rate of interest on silver (first and second) or on barley (third) is 20% per year.

(1) YBC 4669, reverse no. 11

$(1 + 0;12)^{3x} = 1$

In the text the principal of silver (x) is given without its solution, x = 0;34,43,20.

(2) VAT 8528, no. 1

$2^{\frac{x}{5}} = 1,4$

Here the interest is compounded every five years. If the scribe of this tablet had used a table for the logarithmic function $y = \log_2 x$ like MLC 2078 in which

1,4-e 6 íb-si₈ "1,4 corresponds to 6" or "$\log_2 1,4 = 6$"

is listed,[10] he would have obtained the answer more easily:

x/5 = $\log_2 1,4$ = 6, therefore x = 30 (years).

For the clumsy solution performed by the scribe, see my English paper cited in note 9.

(3) AO 6770, no. 2

$(1 + 0;12)^x = 2$

Although the correct answer is

$$x = \frac{\log 2}{\log 6 - \log 5} (\approx 3;48,6),$$

the Babylonian scribe obtained an approximate value of it using linear interpolation:

x ≈ 4 (years) − 2;33,20 (months), that is,

x ≈ 3;47,13,20 (years).

For the details, see my papers cited in note 9.

## §4. Conclusion

We have confirmed that the origin of compound interest goes back 4400 years to Sumer and the custom of compound interest was handed down to the Babylonians. It may be a surprising fact for us that at this early time the Sumerians had already invented what shows both human beings' rational mind and greediness, that is, compound interest. As to its negative aspect, the Babylonians warned of borrowing money easily by their proverbs written in Sumerian, for example.[11]

 ama ír-še$_8$-še$_8$ [--- --- ---] ur$_5$-ra bí-íb-kú-en

 "A weeping mother [said to her son (?)] 'You are swallowed up in debt.'"

Finally, in honour of the late Professor Kramer I would like to conclude my discussion with these words "Compound interest began at Sumer".

## Notes

(1) CDLI-Found; Texts; RIME 1.090501, ex. 01 (http://cdli.ucla.edu/search/search_results. php ? Search Mode = Text & Object ID = P 222532).

 The translation of the text is mine.

(2) Enmetena Foundation Cone A II, 21–26.

(3) R. Jestin, *Tablettes sumériennes de Šuruppak*, 1937.

(4) D. O. Edzard, *Sumerian Grammar*, 2003, pp. 61–66.

(5) O. Neugebauer, *Mathematische Keilschrift-Texte* (MKT), I, II, 1935, III, 1937.

 See BM 85196, no. 9 in MKT II, p. 44.

(6) The text is a fragment of the Code of Hammurabi:

 A. Poebel, *Historical and Grammatical Texts*, 1914, no. 93, obverse column 1, lines

## Acknowledgments

I am very grateful to Prof. K. Maekawa for his kind reply to my inquiry about the books and articles concerning the Enmetena Foundation Cone. I shall, of course, take full responsibility for my interpretation of the passages discussed above.